\documentclass[12 pt]{amsart}
\usepackage{amscd,amssymb}
\usepackage{array}
\oddsidemargin=0.3in \evensidemargin=0.3in \theoremstyle{plain}
\newtheorem{thm}[subsection]{Theorem}

\newtheorem{prop}[subsection]{Proposition}
\newtheorem{cor}[subsection]{Corollary}
\theoremstyle{definition} \theoremstyle{procedure}
\newtheorem{rem}[subsection]{Remark}
\newtheorem{defn}[subsection]{Definition}
\newtheorem{exmp}[subsection]{Example}
\newtheorem{alg}{Algorithm}

\begin{document}
\title{Standard bases over rings}
\author{Afshan Sadiq$^*$}
\address{*Abdus Salam School of Mathematical Sciences, GC University, Lahore, Pakistan}
\email{afshanatiq@gmail.com}

\dedicatory{Dedicated to my Fianc\'{e} Muhammad Atiq Jamil }
\maketitle
\begin{abstract}
The theory of standard bases in polynomial rings with coefficients in a ring $R$ with respect to local orderings is developed. $R$ is a commutative Noetherian ring with $1$ and we assume that linear equations are solvable in $R$.
\end{abstract}
\section{Introduction}
The aim of this paper is to develop the theory of standard bases especially for non-global orderings for polynomial rings with coefficients in a ring. We generalize the concept of Adams and Loustaunau $([1])$ and Greuel and Pfister $([2])$. In the book of Adams and Loustaunau the concept of Gr\"{o}bner bases over polynomial rings with coefficients in a ring is developed, i.e, they consider standard bases with respect to global orderings. In the book of Greuel and Pfister the concept of standard bases over polynomial rings with coefficients in a field is developed, i.e, they consider also non-global orderings. We will generalize both concepts to a uniform theory. \\Note that the theory of standard bases for ideals developed in this paper can also be also developed for modules without changing the proofs. \\First of all we will prove that in the general case the computation of a standard basis with respect to a non-global ordering can be reduced using homogenization to the computation of a Gr\"{o}bner basis with respect to a suitable global ordering. This is also here a very expensive way to compute a standard basis. Therefore later a more efficient algorithm similar to $[2]$ is presented. \\Standard basis computations over the rings $\mathbb{Z}$ and $\mathbb{Z}/\!\!<\!m\!>\!$ can be performed using the computer algebra system SINGULAR (cf. [4]). \\Standard bases are useful in computing elimination of variables, intersection of ideals, quotient of ideals, kernel of the ring map. \\This can be done using the method described in $[2]$, all the results remain the same in our case.

\section{Basic Definitions}
\noindent Let $R$ be a Noetherian commutative ring with $1$ and $R[x_1, \ldots ,x_n]$ the polynomial ring in $n$ variables with coefficients in $R$. Assume that linear equations are $solvable$ in $R$.
\begin{defn}
Linear equations are $solvable$ in $R$ if the following conditions are satisfied for any $a, a_1, \ldots ,a_m \in R$:
\\(1) there is an algorithm to compute generators for the $R$-module $$syz_R(a_1, \ldots ,a_m)=\{(b_1, \ldots ,b_m)\in R^m\; |\;a_1b_1+\ldots +a_mb_m=0\},$$
\\(2) there is an algorithm to determine if $a\in \langle a_1, \ldots ,a_m\rangle$,
\\(3) there is an  algorithm to compute $b_1, \ldots ,b_m \in R$ such that $a=b_1a_1+\ldots +b_ma_m$ if $a\in \langle a_1,\ldots ,a_m\rangle $.
\end{defn}
We will use the notations from $[1]$ and $[2]$ and repeat them here for the convenience of the reader.
\begin{defn}
A $monomial\,ordering$ $>$ is a total ordering on the set of monomials $Mon_n\!=\!\{x^\alpha \!=x_1^{\alpha_1}\cdot \ldots \cdot x_n^{\alpha_n} | \alpha =(\alpha_1,\ldots ,\alpha_n) \in \mathbb{N}^n\}$ in $n$ variables satisfying $$x^\alpha > x^\beta \Longrightarrow x^\gamma x^\alpha > x^\gamma x^\beta$$
for all $\alpha, \beta, \gamma  \in \mathbb{N}^n$. We also say, $>$ is a monomial ordering on $R[x_1,\ldots, x_n]$, if $>$ is a monomial ordering on $Mon_n$.
\begin{exmp}
The local lexicographical ordering $>_{ls}$ on $\operatorname{Mon}_n$ is defined as follows,
\\$x^\alpha > x^\beta$ $\Longleftrightarrow$ $\exists \; 1\leq i \leq n$, $\alpha_1=\beta_1, \ldots ,\alpha_{i-1}=\beta_{i-1},\; \alpha_i < \beta_i$.
\end{exmp}
\end{defn}
\begin{exmp}
Let $M$ be an invertible $(n\times n)$-matrix with real coefficients and $M_1,\ldots ,M_n$ the rows of $M$. The matrix $M$ defines an ordering $>$ on $\operatorname{Mon}_n$ as follows:
\\ \indent \quad  $x^\alpha > x^\beta$ $\Longleftrightarrow$ $\exists \; 1\leq i \leq n$, $M_1\alpha=M_1\beta, \ldots ,M_{i-1}\alpha=M_{i-1}\beta,\; M_i\alpha > M_i\beta$. \\Every ordering can be defined by a matrix $(cf.\, [2])$.
\end{exmp}
\begin{defn}
Let $>$ be a fixed monomial ordering.  Writing $f\in R[x_1, \ldots , x_n]$, $f\neq 0$, in a unique way as a sum of non-zero terms
$$f=a_{\alpha_1} x^{\alpha_1} +a_{\alpha_2} x^{\alpha_2} +\ldots +a_{\alpha_s} x^{\alpha_s} , \, \,\,\,\,\,\,\,\,\, x^{\alpha_1} > x^{\alpha_2} > \ldots  > x^{\alpha_s},$$ and $a_{\alpha_1},a_{\alpha_2} ,\ldots ,a_{\alpha_s} \in R$. We call:
\\$LM(f):=x^{\alpha_1}$, the $leading$ $monomial$ of $f$,
\\$LE(f):={\alpha_1}$, the $leading$  $exponent$ of $f$,
\\$LT(f):=a_{\alpha_1} x^{\alpha_1}$, the $leading$ $term$ of $f$,
\\$LC(f):=a_{\alpha_1}$, the $leading$ $coefficient$ of $f$,
\\$ecart(f):=deg(f)-deg(LM(f))$.
\\We define the leading monomial and the leading term of $0$ to be $0$, and $0$ to be smaller than any monomial.

\end{defn}
\begin{defn}
Let $>$ be a monomial ordering on $\operatorname{Mon}_n$,
$>$ is a called global $($resp. local$)$ ordering if $x^\alpha > 1$ $($resp. $x^\alpha< 1)$ for all $\alpha \neq (0, \ldots ,0)$.

\end{defn}
\begin{defn}
For any monomial ordering $>$ on $\operatorname{Mon}_n$, $$S_>:=\{u\in R[x_1, \ldots , x_n]\backslash \{0\}|LT(u)=1 \}$$
is a $multiplicatively \, closed \, set$, $$R[x_1, \ldots , x_n]_>:=S_>^{-1}R[x_1, \ldots , x_n]=\{\frac{f}{g}|f,g \in R[x_1, \ldots , x_n] \, and \, g\in S_>\}$$ is the localisation of $R[x_1, \ldots , x_n]$ with respect to $S_>$ and we call $R[x_1, \ldots , x_n]_>$ the ring associated to $R[x_1, \ldots , x_n]$ and $>$.
\end{defn}
\begin{defn}
Let $>$ be any monomial ordering.
\\ For $f \in R[x_1, \ldots , x_n]_>$, choose $u \in R[x_1, \ldots , x_n]$ such that $LT(u)=1$ and $uf \in R[x_1, \ldots , x_n]$. Then $$LM(f):=LM(uf),$$
$$LC(f):=LC(uf),$$
$$LT(f):=LT(uf),$$
$$LE(f):=LE(uf).$$
\begin{defn}
Let $>$ be any monomial ordering then, for each $G \subset R[x_1, \ldots , x_n]_>$ $$L(G)=\langle \{LT(g)\,| \, g\in G\}\rangle_{R[x_1, \ldots , x_n]}$$ is called the leading ideal of $G$.
\end{defn}
\end{defn}
\begin{defn}
Let $I\;< \; R[x_1, \ldots , x_n]_>$. \\(1). A finite set $G \subset R[x_1, \ldots , x_n]_>$ is called a standard basis of $I$ with respect to $>$ if $G \subset I$, and $L(I)=L(G)$. \\(2). $G$ is called a strong standard basis\footnote{Strong standard bases do not exist in general. They exist always if $R$ is a principal ideal domain, (cf. theorem 6.4). A strong standard basis is a standard basis. }, of $I$ with respect to $>$, if $G\subset I$ and for any $f\in I\backslash \{0\}$ there exists $i\in \{1, \ldots, t\}$ such that $LT(g_i)$ divides $LT(f)$. \\(3). If $>$ is global, a standard basis is also called a Gr\"{o}bner basis.
\end{defn}

\section{Computing Standard Bases By Using Homogenization}
\begin{thm}
Let $f_1,\ldots ,f_m \in R[x_1,\ldots,x_n]$ generating the ideal $I<R[x_1,\ldots,x_n]_>$, where $>$ is a monomial ordering given by a matrix $M$. Let $F_i:=f_i ^h \in R[t,x_1,\ldots,x_n]$ be the homogenization of $f_i$ and  $>_h$ be the monomial ordering given by the matrix
$\left(
  \setlength{\extrarowheight}{-5pt}\begin{array}{cccc}
    1 & 1 & \cdots & 1 \\
    0 &  &  &  \\
    \vdots &  & M  &  \\
    0 &  &  &  \\
  \end{array}
\right)$. \\Let $\{G_1,\ldots ,G_s \}$ be a Gr\"{o}bner basis, respectively strong Gr\"{o}bner basis of \\$J\!=\langle F_1, \ldots ,F_m\rangle$ with respect to $>_h$. If we denote $g_i=G_i|_{t=1}$ then, $\{g_1, \ldots ,g_s\}$ is a standard basis, respectively strong standard basis of the ideal $I$ with respect to $>$.
\end{thm}
\begin{proof}
Assume that $\{G_1,\ldots ,G_s \}$ is a $Gr\ddot{o}bner\, basis$ with respect to $>_h$. \\Let $f\in I\cap R[x_1,\ldots,x_n]$. Then there exists $u\in S_>$ and $\eta _i\in R[x_1,\ldots,x_n] $ such that
$$u\cdot f=\sum_{i=1}^m\eta _i\cdot f_i.$$
Then there exists $\rho,\; \rho_i \in \mathbb{Z}, \rho,\;\rho_i \geq 0$ such that $$t^{\rho}\cdot u^{h}\cdot f^{h}=\sum _{i=1}^m t^{\rho_i}\cdot \eta ^{h} _i\cdot f_i^h=\sum _{i=1}^m t^{\rho_i} \cdot \eta ^{h} _i\cdot F_i.$$
As $\{G_1,\ldots ,G_s \}$ is a Gr\"{o}bner basis $J$ so there exist $\xi_i\in R[x_1,\ldots,x_n]$ such that $$LT(t^{\rho}\cdot u^{h}\cdot f^{h})=\sum_{i=1}^s\xi_i \cdot LT(G_i)$$
putting $t=1$ we obtain the result. \\Now assume that $\{G_1,\ldots ,G_s \}\subseteq J$ is a strong Gr\"{o}bner basis with respect to $>$ and let $f\in I\cap R[x_1,\ldots,x_n]$. We want to show there exists $i$ such that $LT(g_i)|LT(f)$ and that $g_1,\ldots ,g_s\in IR[x_1,\ldots,x_n]_>$. \\As $\{G_1,\ldots ,G_s \}\subseteq J$  $$G_i=\sum_{j=1}^m \xi _{i,j}\cdot F_j$$
with $\xi _{i,j} \in R[t,x_1,\ldots,x_n]$. \\Put $t=1$, we get $$g_i=\sum_{j=1}^m \xi _{i,j|_{t=1}}\cdot f_j$$
this implies $g_1,\ldots ,g_s\in I R[x_1,\ldots,x_n]_>$. \\Now for $f$ there exists $w\in S_>$ such that $$w\cdot f=\sum_{j=1}^m \eta _j\cdot f_j $$
for suitable $\eta_j \in R[x_1,\ldots,x_n]$. Then there exists $\rho,\;\rho_i \in \mathbb{Z},\rho,\;\rho_i \geq 0$ such that $$t^{\rho}\cdot w^{h}\cdot f^{h}=\sum t^{\rho_i}\cdot \eta ^{h} _i\cdot f^{h} _j=\sum _{j=1}^m t^{\rho_i}\cdot \eta ^{h} _i\cdot F_j.$$
As $\{G_1,\ldots ,G_s \}$ is a strong Gr\"{o}bner basis of $J$ there exists $i$ such that \\$LT(G_i)|LT(t^{\rho}\cdot w^{h}\cdot f^{h})$. This implies $LT(g_i)|LT(f)$.
\end{proof}
\section{Normal Form}

The concept of a normal form with respect to a given system of polynomials is the basis of the theory of standard bases. Normal forms for non-global orderings are different and more complicated than normal forms for global orderings. This is already the case for polynomial rings over a field.
\begin{defn}
Let $\mathcal{G}$ denote the set of all finite lists $G\subset R[x_1, \ldots , x_n]_>$,
$$NF:R[x_1, \ldots , x_n]_>\times \mathcal{G} \longrightarrow R[x_1, \ldots , x_n]_>, \, \, (f,G)\longmapsto NF(f|G),$$
is called a $ normal \, form$ on $R[x_1, \ldots , x_n]_>$ if for all $G\subset \mathcal{G}$ and for all \\$f\in R[x_1, \ldots , x_n]_> $,
\\(1). $NF(0|G)=0$,
\\(2). $NF(f|G) \neq 0 \, \Longrightarrow LT(NF(f|G)) \notin L(G).$
\\(3). If $G=\{g_1, \ldots , g_s\}$, then there exists $u\in S_>$ such that $r:=uf-NF(f|G)$ \\has a $standard \, representation$ with respect to $G$, that is, $$r=\sum _{i=1}^s \xi _i\cdot g_i$$
for suitable $\xi_i \in R[x_1, \ldots , x_n]$ and $LM(r)=max_{i=1}^s \{LM( \xi _i)LM( g_i)\}$.
\end{defn}
To prove the existence of a normal form we give an algorithm to compute it.
\begin{defn}
Let $h\in R[x_1, \ldots , x_n]$, $T\subseteq R[x_1, \ldots , x_n]$ be finite. If $h=0$, let $S(T,h)=\emptyset$. If $h\neq 0$, let \\$S(T,h):=\{\sum _{g\in T}c_g\cdot x^{\alpha_g}\cdot g |\;LT(\sum _{g\in T}c_g\cdot x^{\alpha_g}\cdot g)=LT(h),\, c_g\in R $ and $x^{\alpha_g }\cdot LM(c_g\cdot g)=LM(h)$ if $c_g\neq0\}$.
\end{defn}
\begin{rem}
The set $S(T,h)$ can be infinite. Algorithm $1$ requires to choose an element of $S(T,h)$ which is of minimal ecart. This is achieved by computing a generating system of $S(T,h)$, which is a kind of a syzygy module, so by the assumption on $R$ this can be done.
\end{rem}
\begin{alg}
$NF(f|G)$\\Let $>$ be any monomial ordering.
\texttt{\\Input}: $f\in R[x_1, \ldots , x_n]$, $G=\{g_1,\ldots ,g_s\}\subset R[x_1, \ldots , x_n]$ with $g_i\neq 0 \;\;\forall \;i=1,\ldots ,s$.
\texttt{\\Output}: $h\in R[x_1, \ldots , x_n]$ a normal form of $f$ with respect to $G$.\\$\bullet$ $h:=f$;
\\$\bullet$ $T:=G$;
\\$\bullet$ while$(S(T,h)\neq \emptyset )$
\\ \indent \quad choose $k\in S(T,h)$ such that $ecart(k)$ is minimal;
\\ \indent \quad if$(ecart(k)>ecart(h))$
\\ \indent \qquad $T:=T \cup \{h\}$;
\\ \indent \quad $h:=h-k$;
\\ $\bullet$ return $h$;
\end{alg}
\begin{prop}
The algorithm terminates and defines a normal form.
\end{prop}
\begin{proof}
Termination is proved by using homogenization with respect to $t$: \\We start with $h:=f^h$ and $T^h:=\{g^h|g\in G\}$. \\The while loop looks as follows
\\$\bullet$ while$(S(T^h,t^\alpha h)\neq \emptyset$ for some $\alpha)$
\\ \indent \quad choose $k\in S(T^h,t^\alpha h)$ such that $\alpha \geq 0$ is minimal;
\\ \indent \quad if$(\alpha > 0)$
\\ \indent \qquad $T^h:=T^h \cup \{h\}$;
\\ \indent \quad $h:=t^{\alpha}h-k$;
\\  \indent \quad $h:=(h|_{t=1})^h$;
\\ Since $R[x_1, \ldots , x_n]_>$ is Noetherian, there exists some positive integer $N$ such that $L(T^h_v)$ becomes stable for $v \geq N$, where $T^h_v$ denotes the set $T^h$ after the $v$-th turn of the while loop. The next $h$, therefore, satisfies $LT(h)\in L(T^h_N)=L(T^h)$, whence, $LT(\sum _{g\in T^h}c_g x^{\alpha_g}g )\\=LT(h)$ for some $\sum _{g\in T^h}c_g x^{\alpha_g}g$ and $\alpha =0$. That is, $T^h_v$ itself becomes stable for $v \geq N$ and the algorithm continues with fixed $T^h$. Then it terminates, since $>_h$ is a well ordering on $R[t,x_1, \ldots , x_n]$.
For the correctness consider the $v$-th while loop of Algorithm $1$. There we create (with $h_0:=f$) $$h_v=h_{v-1}-\sum _{g\in T}c_g x^{\alpha _g} g$$ for some $\sum _{g\in T}c_g x^{\alpha _g} g $ such that $LT(\sum _{g\in T} x^{\alpha_g}g )=LT(h_{v-1})$ and $x^{\alpha_g}LM(c_gg)=LM(h_{v-1})$ if $c_g\neq 0$, from the construction of $T^h$ we have $$\sum _{g\in T}c_g x^{\alpha_g}g = \sum_{i=1}^sc_ix^{\alpha_i}g_i+\sum_{j=0}^{v-2}d_jx^{\beta_j}h_j$$
where $$c_g=\left\{
    \begin{array}{ll}
      c_i & \,\, g=g_i\\
      d_j & \,\, g=h_j
    \end{array}
  \right.
$$
$$\alpha_g=\left\{
    \begin{array}{ll}
      \alpha_i & \,\, g=g_i \\
      \beta_j & \,\, g=h_j
    \end{array}
  \right.
$$
which implies
$$h_v=h_{v-1}-(\sum_{i=1}^sc_ix^{\alpha_i}g_i+\sum_{j=0}^{v-2}d_jx^{\beta_j}h_j).$$
Especially for $v\geq 2$ we have $LM(f)>LM(h_{v-1})=LM(x^{\alpha_i}c_ig_i)$ if $c_i\neq 0$ and $LM(f)>LM(h_{v-1})$ \\$=LM(x^{\beta_j}d_jh_j)$ if $d_j\neq 0$. This implies especially $x^{\beta_j}<1$. \\Suppose by induction, that in the first $v-1$ steps $(v\geq 1)$ we have constructed standard representations $$u_jf=\sum_{i=1}^sa_i^{(j)}g_i+h_j, \, \, u_j\in S_>, \, \, a_i^{(j)}\in R[x_1, \ldots , x_n],$$
with $LM(u_jf-h_j)=max_{i=1}^s  \{LM(a_i^{(j)})LM(g_i)\}$ for $0\leq j \leq v-1$, starting with $u_0=1, \, a_i^{(0)}=0$. \\Consider this standard representation for $j=v-1$. \\We replace $h_{v-1}$ by $h_v+(\sum_{i=1}^sc_ix^{\alpha_i}g_i+\sum_{j=0}^{v-2}d_jx^{\beta_j}h_j)$, hence we obtain
$$u_{v-1}f=\sum_{i=1}^sa_i^{(v-1)}g_i+h_v+(\sum_{i=1}^sc_ix^{\alpha_i}g_i+\sum_{j=0}^{v-2}d_jx^{\beta_j}h_j),$$
where each $h_j$ has a standard representation as above
$$u_{v-1}f=\sum_{i=1}^sa_i^{(v-1)}g_i+h_v+(\sum_{i=1}^sc_ix^{\alpha_i}g_i+\sum_{j=0}^{v-2}d_jx^{\beta_j}(u_jf-\sum_{i=1}^sa_i^{(j)}g_i))$$
$$(u_{v-1}-\sum_{j=0}^{v-2}d_jx^{\beta_j}u_j)f=
(\sum_{i=1}^sa_i^{(v-1)}g_i+\sum_{i=1}^sc_ix^{\alpha_i}g_i-\sum_{j=0}^{v-2}d_jx^{\beta_j}\sum_{i=1}^sa_i^{(j)}g_i)+h_v.$$
$$=\sum_{i=1}^s(a_i^{(v-1)}+c_ix^{\alpha_i}-\sum_{j=0}^{v-2}d_jx^{\beta_j}a_i^{(j)})g_i+h_v.$$
Let $u_v:=(u_{v-1}-\sum_{j=0}^{v-2}d_jx^{\beta_j}u_j)$ and $a_i^{(v)}:=a_i^{(v-1)}+c_ix^{\alpha_i}-\sum_{j=0}^{v-2}d_jx^{\beta_j}a_i^{(j)}$. We have to prove that $u_v\in S_>$ and $$u_vf=\sum_{i=1}^sa_i^{(v)}g_i+h_v$$
is a standard representation, i.e, \\$LM(f)=LM(u_vf-h_v)=\max_{i=1}^s\{LM(a_i^{(v)})LM(g_i)\}$. Since $x^{\beta_j}<1$ in case $d_j\neq 0$ it follows $u_v\in S_>$.
\\Since $LM(f)>LM(x^{\alpha_k}c_kg_k)$ if $c_k\neq 0$ and $LM(f)\geq LM(a_k^{(j)}g_k)$, $x^{\beta_j}<1$, it follows $LM(a_k^{(v)}g_k)\leq LM(f)$. If $LM(a_i^{(v-1)}g_i)=LM(f)$ then with the same argument we obtain $LM(a_i^{(v)}g_i)=LM(a_i^{(v-1)}g_i)=LM(f)$.
\end{proof}
\begin{exmp}
We consider $R=\mathbb{Z}$ and use the local lexicographical ordering $ls$ with $x>y$ in $\mathbb{Z}[x,y]$. Let $f=xy^4-12x^2$  then $ecart(f)=0$ and let $G=\{f_1,f_2,f_3\}$ where $$f_1=-3x+xy,f_2=y^2-2x^2y,\,f_3=6x^2-x^3y^2$$
then $ecart(f_1)=1,ecart(f_2)=1,ecart(f_3)=3$. \\In step 1:
\\$h_0=xy^4-12x^2$, $T:=G$ and $$S(T,h_0)=\{ ky^4f_1+(3k+1)xy^2f_2, \,\,k\in \mathbb{Z}\}.$$ All elements in $S(T,h_0)$ have ecart $1$ and we choose $xy^4-xy^5+4x^3y^3\in S(T,h_0)$. Since $ecart(h_0)<ecart(xy^4-xy^5+4x^3y^3)$ we have to enlarge $T$: \\ $T=T\cup \{f_4:=h_0\}$ and \\$h_1=xy^4-12x^2-(xy^4-xy^5+4x^3y^3)=xy^5-12x^2-4x^3y^3$. \\In step 2: \\$h_1=xy^5-12x^2-4x^3y^3$ with $ecart(h_1)=0$, $$S(T,h_1)=\{ky^5f_1+lxy^3f_2+(1+3k-l)yf_4,\,\, k,\,l\in \mathbb{Z}\}.$$
We choose $xy^5-12x^2y\in S(T,h_1)$ with minimal ecart $0$ and obtain \\$h_2=xy^5-12x^2-4x^3y^3-(xy^5-12x^2y)=-12x^2+12x^2y-4x^3y^3$. \\In step 3: \\$h_2=-12x^2+12x^2y-4x^3y^3$ with $ecart(h_2)=4$, $$S(T,h_2)=\{(2k+4)xf_1+kf_2, \,\, k\in \mathbb{Z}\}.$$
We choose $-12x^2+4x^2y\in S(T,h_2)$ with minimal ecart $1$ and obtain \\$h_3=-12x^2+12x^2y-4x^3y^3-(-12x^2+4x^2y)=8x^2y-4x^3y^3$. \\In step 4: \\$h_3=8x^2y-4x^3y^3$ and $LT(h_3)\notin L(T)$, thus $NF(f|G)=8x^2y-4x^3y^3$.
\end{exmp}
\begin{rem}
Assume $R$ has the following property: $c=a_1x_1+\ldots +a_sx_s$ is solvable in $R$ if and only if there exists $j$ and $x\in R$ such that $c=a_jx$. \\Then normal form algorithm is similar to the corresponding normal form algorithm for a polynomial ring over a field, i.e, $S(T,h)$ can be replaced by $S(T,h)=\{g\in T|\;LT(g)|LT(h)\}$. In this case each standard basis is a strong standard basis. If $R$ is a discrete valuation ring or $R=\mathbb{Z}/\langle p^n\rangle$, $p$ a prime number, then $R$ has the property above.
\end{rem}

\section{Computing Standard Bases}

\begin{thm}
Let $I$ $<$ $R[x_1, \ldots , x_n]_>$ and let $G=\{g_1 , \ldots , g_t\}$ be a set of non-zero polynomials in $I$. Then the following are equivalent.
\\(1). $L(G)=L(I)$.
\\(2). For any polynomial $f\in R[x_1, \ldots , x_n]_>$, we have $f\in I$ if and only if \\ \indent \quad $NF(f|G)\!=\!0$.
\\(3). For all $f\in I, \, uf=\sum_{i=1}^th_ig_i$ for some polynomials $u,h_1, \ldots , h_t \in R[x_1, \ldots , x_n]$ \\ \indent \quad such that
$LT(u)=1$ and \\ \indent \quad $LM(f)=max_{i=1}^t \{LM( h_i)LM( g_i)\}$.
\end{thm}
\begin{proof}
$(1) \Longrightarrow (2).$ We know that if $NF(f|G)=0$, then $f\in I$. Conversely assume that $f\in I$. Let $r=NF(f|G)$ and assume $r \neq 0$. Since $G \subset I$ we have $r\in I$. This implies $LT(r) \in L(I)=L(G)$. This is a contradiction to Definition $9$.
\\$(2) \Longrightarrow (3).$ This is obvious from the Definition of normal form.
\\ $(3) \Longrightarrow (1).$ For $f\in I$ we need to show that $LT(f)\in L(G)$. We have that $uf=\sum_{i=1}^th_ig_i$ such that $LM(f)=max_{i=1}^t \{LM(h_i)LM(g_i)\}$. It is easily seen that $LT(f)=\sum LT(h_i)LT(g_i)$ where the sum is over all $i$ satisfying $LM(f)=LM(h_i)LM(g_i)$. Thus $LT(f)\in L(G)$, as desired.
\end{proof}

\begin{cor}
If $G$ is a standard basis of $I$ then $I$ is generated by $G$.
\end{cor}
\begin{proof}
Clearly $\langle g_1,\ldots, g_t\rangle \subset I$, since each $g_i$ is in $I$. For the other inclusion, let $f\in I$.  By Theorem $5.1$, $NF(f|G)=0$ and hence $uf\in \langle g_1, \ldots ,g_t \rangle_{R[x_1, \ldots , x_n]}$ for a suitable $u\in S_>$.
\end{proof}
\begin{rem}
Let $K$ be a field and $R=K[y_1,\ldots ,y_m]$ be the polynomial ring with variables $y_1,\ldots,y_m$. Let $>$ be a product ordering on $K[y_1,\ldots,y_m,x_1,\ldots,x_n]$ such that the $x_1,\ldots,x_n$ dominate $y_1,\ldots,y_m$ and the restriction of $>$ to $R$ is global. Let $I< K[y_1,\ldots,y_m,x_1,\ldots,x_n]_>$ an ideal and $G=\{f_1,\ldots,f_s\}$ a standard basis of $I$. \\Let $>_1$ be the ordering on $R[x_1,\ldots,x_n]$ considered as polynomial ring with coefficients in $R$ induced by $>$. Then $G$ is a standard basis of $I$ with respect to $>_1$.
\end{rem}
\begin{proof}
Let $f\in I$ then $NF_>(f|G)=0$ ($NF_>$ the normal form with respect to $>$). Analyzing the algorithm for $NF_>$ we obtain $a_1,\ldots,a_s\in R$ such that $\sum_{i=1}^sa_iLT_{>_1}(f_i)=LT_{>_1}(f)$. Here $LT_{>_1}$ is the leading term with respect to the ordering $>_1$.
\end{proof}

We use the following Definition from $[1]$.
\begin{defn}
Given monomials $x^{\alpha _1}, \ldots, x^{\alpha _s}$ and non-zero elements $c_1, \ldots, c_s$ in $R$ set $L=(c_1x^{\alpha _1}, \ldots, c_sx^{\alpha _s})$. Then for a given monomial $x^\alpha$, we call a $syzygy$ $h=(h_1,\ldots ,h_s)\in syz(L) \subset (R[x_1, \ldots, x_n])^s$ homogeneous of degree $x^{\alpha}$ provided that each $h_i$ is a term and $x^{\alpha_i}LM(h_i)=x^{\alpha}$ for all $i$ such that $h_i \neq 0$.
\end{defn}
\begin{thm}
Let $G=\{g_1, \ldots , g_t\}$ be a set of non-zero polynomials in $R[x_1, \ldots , x_n]$. Let $\mathcal{B}$ be a homogeneous generating set for
$syz(LT(g_1), \ldots , LT(g_t))$. Then $G$ is a standard basis for the ideal $\langle g_1, \ldots , g_t\rangle R[x_1, \ldots , x_n]_>$ if and only if for all \\$(h_1, \ldots , h_t)\in \mathcal{B}$, we have $$NF(\sum_{i=1}^th_ig_i|G)=0.$$

\end{thm}
\begin{proof}
If $G$ is a standard basis, then by Theorem $5.1$, $$NF(\sum_{i=1}^t h_ig_i|G)=0.$$
\\Conversely, let $g\in \langle g_1, \ldots ,g_t\rangle R[x_1, \ldots , x_n]_>$, then there exists $u\in S_>$ and $u_i\in R[x_1, \ldots , x_n]$ such that $$ug=\sum_{i=1}^t u_ig_i.$$
Choose a representation as in above equation with $x^\alpha = max_{i=1}^t(LM(u_i)LM(g_i))$ minimal. Since by Theorem $5.1$, we need to show that $LM(g)=x^\alpha$. We assume $LM(g)<x^\alpha$ and show that we can obtain an equation for $g$ with a smaller value for $x^\alpha$. Let $S=\{i\in\{1, \ldots, t\}|LM(u_i)LM(g_i)=x^\alpha\}$. Then $$\sum_{i\in S}LT(u_i)LT(g_i)=0.$$
Let $\textbf{h}=\sum_{i\in S}LT(u_i)\textbf{e}_i$ (where $\textbf{e}_1=(1,\ldots ,0), \ldots ,\textbf{e}_t=(0,\ldots ,1)$ is a generating set for $R[x_1, \ldots , x_n]^t$). Then $\textbf{h}\in syz(LT(g_1),\ldots ,LT(g_t))$ and $\textbf{h}$ is homogeneous of degree $x^\alpha$. Now let $\mathcal{B}=\{\textbf{h}_1, \ldots ,\textbf{h}_l\}$, with $\textbf{h}_j=(h_{1,j}, \ldots ,h_{t,j})$ then $\textbf{h}=\sum_{j=1}^l a_j \textbf{h}_j$. Since $\textbf{h}$ is a homogeneous syzygy, we may assume that the $a_j$ are terms such that $LM(a_j)LM(h_{i,j})LM(g_i)=x^\alpha$ for all $i,j$ such that $a_jh_{i,j}\neq 0$. By hypothesis, for each $j, \, NF(\sum_{i=1}^th_{i,j}g_i|G)=0$. Thus by Theorem $5.1$, for each $j=1, \ldots ,l$ there exist $w_j\in S_>$ and $v_{i,j}\in R[x_1, \ldots , x_n]$ such that $$w_j\sum_{i\in S}h_{i,j}g_i=\sum_{i=1}^tv_{i,j}g_i,$$
and $$max_{i=1}^tLM(v_{i,j}g_i)=LM(\sum_{i=1}^th_{i,j}g_i)<max_{i=1}^tLM(h_{i,j})LM(g_i).$$
\\The latter strict inequality is because $\sum_{i=1}^th_{i,j}LT(g_i)=0.$ \\We may assume $w=w_j$ for all $j$. \\Thus, $$wug=\sum_{i=1}^t wu_ig_i$$
$$=\sum_{i\in S} wLT(u_i)g_i+\sum_{i\in S}(wu_i-wLT(u_i))g_i+\sum_{i\notin S}wu_ig_i$$
$$=\sum_{j=1}^l\sum_{i\in S}wa_jh_{i,j}g_i \,+\, terms \,lower\, than\, x^\alpha$$
$$=\sum_{j=1}^l\sum_{i=1}^ta_jv_{i,j}g_i \,+\, terms\, lower\, than \,x^\alpha.$$
We have $max_{i,j}LM(a_j)LM(v_{i,j})LM(g_i)<max_{i,j}LM(a_j)LM(h_{i,j})LM(g_i)=x^\alpha.$ We have a representation of $g$ as a linear combination of the $g_i$ such that the maximum of the leading monomials of any summand is less than $x^\alpha$. Thus the theorem is proved.
\end{proof}
As a consequence of Theorem $5.5$, we obtain that the following algorithm computes a standard basis for a given ideal $I$ in $R[x_1, \ldots , x_n]_>$.
\begin{alg}Standardbasis$(G)$
\texttt{\\Input}: $F=\{f_1,\ldots ,f_r\} \subset R[x_1,\ldots,x_n]$ with $f_i\neq 0\,\,(1\leq i\leq r)$,
\texttt{\\Output}: $G$ a Standard basis for $\langle f_1,\ldots ,f_r\rangle R[x_1,\ldots,x_n]_>$.
\\ $\bullet$ $G:=F$;
\\ $\bullet$ $P$:= a finite homogeneous generating set (considered as ordered set) for \\ \indent $syz(\{LT(f_i)\}_{1\leq i\leq r})$;
\\ $\bullet$ while$(P\neq \emptyset)$
\\ \indent \quad Let $G=\{f_1,\ldots ,f_k\}$
\\ \indent \quad choose $(s_1,\ldots ,s_k)\in P$
\\ \indent \quad $P:=P\backslash \{(s_1,\ldots ,s_k)\}$;
\\ \indent \quad $h:=NF(\sum_{i=1}^ks_if_i|G)$;
\\ \indent \quad if$(h\neq 0)$;
\\ \indent \qquad $f_{k+1}:=h$;
\\ \indent \qquad $G:=\{f_1,\ldots ,f_{k+1}\}$;
\\ \indent \qquad $H$:= a finite homogeneous generating set for $syz(\{LT(f_i)\}_{1\leq i\leq k+1})$;
\\ \indent \qquad $P:=(P\times \{0\})\cup \{h=(h_1,\ldots ,h_{k+1})\in H|h_{k+1}\neq 0\}$;
\\ $\bullet$ return $G$;
\end{alg}
Using ideas of M.M\"{o}ller $([3])$ Adams and Loustaunau propose a more efficient algorithm for computing Gr\"{o}bner bases. This applies also in our situation with the same proof.
\\We use the following Definition and Theorem from $[1]$.
\begin{defn}
Let $x^{\alpha _1}, \ldots, x^{\alpha _s}$ be a set of monomials. For any subset $J \subseteq \{1, \ldots ,s\}$, set $x^{\gamma _J} =lcm(x^{\alpha _j}|j\in J)$. We say that $J$ is saturated with respect to
$x^{\alpha _1}, \ldots, x^{\alpha _s}$ provided that for all $j\in \{1, \ldots , s\}$ if $x^{\alpha _j}$ divides $x^{\gamma _J}$, then $j\in J$. For any subset $J\in \{1, \ldots ,s\}$ we call the saturation of $J$ the set $J_1$ consisting of all $j\in \{1, \ldots ,s\}$ such that $x^{\alpha _j}$ divides
 $x^{\gamma _J}$. (Note that $x^{\gamma _J}=x^{\gamma _{J_1}}.$)
\end{defn}
\begin{thm}
$($cf. $[1]$, page $214$$)$ Given monomials $x^{\alpha _1}, \ldots, x^{\alpha _s}$ and non-zero elements $c_1, \ldots, c_s$ in $R$. For each set $J \subseteq \{1, \ldots ,s\}$,which is $saturated$ with respect to $x^{\alpha _1}, \ldots, x^{\alpha _s}$, let $\mathcal{B}_J=\{b_{1,J}, \ldots ,b_{\nu_J,J}\}$ be a set of generators of the $R$-module of syzygies $syz_R(c_j|j \in J)$. (Note that each of the vectors $b_{\mu ,J}$ is in the $R$-module $R^{|J|}$, where $|J|$ denotes the cardinality of $J$). For each such $b_{\nu ,J}$, denote its $j$th coordinate, for $j\in J$, by $b_{j\nu ,J}$ and $x^{\gamma _J} =lcm(x^{\alpha _j}|j\in J)$. Set
$$s_{\nu ,J}=\sum_{j\in J} b_{j\nu ,J}\frac{x^{\gamma _J}}{x^{\alpha _j}} e_j$$
$($ Note that each of the vectors $s_{\nu ,J}$ is in $R[x_1, \ldots , x_n]^s$ $)$. Then the vectors
$s_{\nu ,J}$, for $J$ running over all such saturated subsets of $\{1, \ldots ,s\}$, and, $1 \leq \nu \leq \nu_J$, forms a homogeneous generating set for the syzygy module
$syz(c_1x^{\alpha _1}, \ldots, c_sx^{\alpha _s})$.
\end{thm}
\begin{exmp}
We consider $R=\mathbb{Z}$ and let $c_1x^{\alpha_1}=3xy^2,\; c_2x^{\alpha_2}=7xyz,\; c_3x^{\alpha_3}=2y^2z^2$. The saturated subsets of $\{1,2,3\}$ are $\{1\},\{2\},\{3\},\{1,2\}$ and $\{1,2,3\}$. Since $\mathbb{Z}$ is an integral domain, the singletons $\{1\},\{2\},\{3\}$ do not give rise to any non-zero syzygy.
\\For $J=\{1,2\}$ we need to solve in $R=\mathbb{Z}$ the equation $3b_1+7b_2=0$.
The module of all solutions is generated by $(7,-3)$. Since $x^{\gamma_J}=xy^2z$, the corresponding syzygy is \\$s_{\nu ,J}= 7 \frac{xy^2z}{xy^2}$ \textbf{$e_1$}$+ 3 \frac{xy^2z}{xyz}$ \textbf{$e_2$}= $(-7z,3y)$.\\Now for $J=\{1,2,3\}$ we need to solve $3b_1+7b_2+2b_3=0$. The module of all solutions is generated by $(-4,2,-1)$ and $(-7,3,0)$. Then with $x^{\gamma_J}=xy^2z^2$ we obtain the syzygies are \\$s_{\nu ,J}=-4\frac{xy^2z^2}{xy^2}$ \textbf{e$_1$}$+2\frac{xy^2z^2}{xyz}$ \textbf{e$_2$} $-\frac{xy^2z^2}{y^2z^2}$ \textbf{e$_3$}= $(-4z^2,2yz,-x)$. \\$s_{\nu ,J}=-7\frac{xy^2z^2}{xy^2}$ \textbf{e$_1$} + $3\frac{xy^2z^2}{xyz}$ \textbf{e$_2$} = $(-7z^2,3yz,0)$. \\So we obtain that $$syz(3xy^2,7xyz,2y^2z^2)= \langle (-7z,3y,0),(-4z^2,2yz,-x),(-7z^2,3yz,0)\rangle. $$

\end{exmp}

The theorem is the basis of the following modified standard basis algorithm.

\begin{alg}
Standardbasis(G)
\texttt{\\Input}: $F=\{f_1,\ldots ,f_s\}\subseteq R[x_1, \ldots , x_n] $ with $f_i \neq 0$ ($1 \leq i \leq s$),
\texttt{\\Output}: $G$ a standard basis for $\langle f_1,\ldots ,f_s\rangle R[x_1, \ldots , x_n]_>$.
\\$\bullet $ $G:=F$;
\\$\bullet $ $\sigma :=1$;
\\$\bullet $ $m:=s$;
\\$\bullet $ while($\sigma \leq m$)
\\  \indent \quad Compute $S=\{$ subsets of $\{1, \ldots , \sigma\} $, saturated with respect to \\ \indent \qquad $LM(f_1), \ldots , LM(f_s)$, which contain $\sigma $  $\}$;
\\  \indent \quad   for($J\in S$)
\\  \indent \qquad \indent   $x^\gamma :=lcm(LM(f_j)|j\in J)$;
\\  \indent \qquad  \indent  Compute a generating set $\{b_{i,J}, i=1, \ldots ,\mu_J\}$ \\ \indent \qquad \indent \qquad \indent  for $\langle LC(f_j)|j\in J,$ $ j \neq \sigma \rangle_R:\langle LC(f_\sigma )\rangle_R$
\\  \indent \quad  \indent \quad for($i:=1, \ldots ,\mu_J$)
\\  \indent \qquad \indent \quad  \indent \quad Compute $b_j \in R,j\in J,j \neq \sigma$
 \\  \indent \qquad \indent \qquad \indent \qquad such that $\sum  _{j\in J,j \neq \sigma}  b_j LC(f_j) + b_{i,J}LC(f_ \sigma)=0$
\\ \indent \qquad \indent \quad  \indent \quad  $r$:=NF$(\sum _{j\in J,j \neq \sigma} b_j \frac{x^\gamma}{LM(f_j)} f_j + b_{i,J} \frac{x^\gamma}{LM(f_\sigma)} f_\sigma|G)$;
\\ \indent \qquad \indent \quad  \indent \quad if($r \neq 0$)
\\  \indent \qquad \indent \qquad \indent  \qquad \indent    $f_{m+1}:=r$;
\\   \indent \qquad \indent \qquad \indent  \qquad \indent   $G:= G \cup \{f_{m+1}\}$;
\\   \indent \qquad \indent \qquad \indent  \qquad \indent   $m:=m+1$;
\\ \indent \qquad   $\sigma :=\sigma +1$;
\\$\bullet $ return $G$;
\end{alg}
\begin{exmp}
We consider $R=\mathbb{Q}[x,y]$ with the local lexicographical ordering $ls$ in $R[z]$. Let $G=\{f_1,f_2\}$ and $I=\langle G\rangle$ where $$f_1=y-x^3z^2 \,\,and \,\,f_2=x^2y-z$$ Then $ecart(f_1)=2$ and $ecart(f_2)=1.$
\\ In step 1 $(\sigma =1):$
\\ $S=\{\{1\}\}=$ saturated subsets of $\{1\}$ containing $1$.
\\ Since $R$ is a domain and $f_1,\,f_2$ are irreducible we have no non trivial syzygy.
\\ In step 2 $(\sigma =2):$
\\ $S=\{\{1,2\}\}=$ saturated subsets of $\{1,2\}$ containing $2$.
\\ \indent $J=\{1,2\}$.
\\ \indent \quad $x^\gamma =lcm(1,1)=1$
\\ \indent \quad A generating set for $\langle y\rangle:\langle x^2y\rangle$ is $\{1\}$.
\\ \indent \qquad The solution of $yb_1+x^2y=0 $ is $b_1=-x^2$.
\\$h=-x^2(y-x^3z^2)+1(x^2y-z)=-z+x^5z^2$ \\which is reduced with respect to $G=\{f_1,f_2\}$. \\ $f_3:=h$,  \\and $G=\{f_1,f_2,f_3\}$.
\\ In step 3 $(\sigma =3):$
\\ $S=\{\{1,2,3\}\}=$ saturated subsets of $\{1,2,3\}$ containing $3$.
\\ \indent $J=\{1,2,3\}$.
\\ \indent \quad $x^\gamma =lcm(1,1,z)=z$.
\\ \indent \quad A generating set for $\langle y,x^2y\rangle:\langle -1\rangle $ is $\{y\}$.
\\ \indent \qquad The solution of $yb_1+x^2yb_2+x^2y=0$ is $b_1=1,b_2=0$. \\$h=z(y-x^3z^2)+y(-z+x^5z^2)=x^5yz^2-x^3z^3$. \\$NF(h|G)=0$. \\So $G=\{f_1,f_2,f_3\}$ is a standard basis.
\end{exmp}
\begin{rem}
Let $R$ be a local ring of the type $K[y_1, \ldots, y_m]_{\langle y_1, \ldots, y_m\rangle}/I$, where
$I< K[x_1, \ldots, x_n]_{\langle x_1, \ldots, x_n\rangle}$. Then we can compute Hilbert-Samuel function of ideals in $R[x_1,\ldots, x_n]_>$ by using the same method described in $[2]$. Carlo Traverso $[5]$, used the Hilbert-Samuel function to speed up the Buchberger Algorithm. The same method can be applied in our case too.
\end{rem}
\section{Standard Bases Over Principal Ideal Domains}

If $R$ is a principal ideal domain (for short PID) then there is a standard basis algorithm similar to the corresponding algorithm for a polynomial ring over a field (cf. $[2]$, page $54$) with the following notion of the s-polynomial.
\begin{defn}
Let $f,g\in R[x_1,\ldots ,x_n]\backslash \{0\}$. $$spoly(f,g):=\frac{lcm(LT(f),LT(g))}{LT(f)}f - \frac{lcm(LT(f),LT(g))}{LT(g)}g.$$

\end{defn}

This is a consequence of $[1]$, Proposition $4.5.3$ and Theorem $5.5$.

\begin{exmp}
We consider $R=\mathbb{Z}$ with the local lexicographical ordering $ls$ with $x>y$ in $\mathbb{Z}[x,y]$. \\Let $I=\langle f_1,f_2\rangle$ where $f_1=-3y+xy$ and $f_2=y^2-2x$. \\Initialization: $G=\{f_1,f_2\}$, $\mathcal{G}=\{\{f_1,f_2\}\}$.
\\ In step 1:
\\$LT(f_1)=-3y,\, LT(f_2)=y^2$, \\$c=lcm(c_1,c_2)=lcm(3,1)=3,\, x^\gamma =lcm(LM(f_1),LM(f_2))=lcm(y,y^2)=y^2$
$$h:=spoly(f_1,f_2)=\frac{3y^2}{-3y} (-3y+xy) -\frac{3y^2}{y^2}(y^2-2x)=6x-xy^2$$
which is reduced with respect to $G$. \\So $G \cup \{f_3:=h\}$ and $\mathcal{G}=\{\{f_1,f_3\},\,\{f_2,f_3\}\}$.
\\ In step 2:
\\$LT(f_1)=-3y, LT(f_3)=6x$, \\$c=lcm(c_1,c_3)=lcm(3,6)=6, \, x^\gamma =lcm(LM(f_1),LM(f_3))=lcm(y,x)=xy$
$$spoly(f_1,f_3)=\frac{6xy}{-3y} (-3y+xy) -\frac{6xy}{6x}(6x-xy^2)=xy^3-2x^2y=xy\cdot f_2$$
$NF(xy^3-2x^2y|G)=0$.
\\ In step 3:
\\$LT(f_2)=y^2, LT(f_3)=6x$, \\$c=lcm(c_2,c_3)=lcm(1,6)=6, \, x^\gamma =lcm(LM(f_2),LM(f_3))=lcm(y^2,x)=xy^2$
$$spoly(f_2,f_3)=\frac{6xy^2}{y^2} (y^2-2x) -\frac{6xy^2}{6x}(6x-xy^2)=xy^4-12x^2$$
$NF(xy^4-12x^2|G)=0$. Since $xy^4-12x^2=xy^2f_2-12x^2+2x^2y^2=xy^2f_2+2xf_3$. \\So $G=\{f_1,f_2,f_3\}$ is a standard basis.
\end{exmp}
\begin{exmp}
We consider $R=\mathbb{Z}$ with the local degree lexicographical ordering $ds$ $(cf.[2])$ with $x>y>z$ in $\mathbb{Z}[x,y,z]$. \\Let $I=\langle f_1,f_2,f_3,f_4\rangle$ where $f_1=15x^2+28y^2z^6$ and $f_2=3x^2y+7yz^5$ $f_3=4xy^2-5xz^{10}$, $f_4=-28y^3+35yz^{11}$.
\\Similar to Example $6.2$ we have to compute the normal form of the spoly's of all pairs $(f_i,f_j)$. \\$Nf(spoly(f_1,f_2)|G)=35yz^5-28y^3z^6=:f_5$, $G=G\cup \{f_5\}$
\\$Nf(spoly(f_2,f_3)|G)=28y^2z^5+15x^2z^{10}=:f_6$, $G=G\cup \{f_6\}$
\\$Nf(spoly(f_3,f_5)|G)=35xz^{15}=:f_7$, $G=G\cup \{f_7\}$.
\\The normal form of all the other spoly's is zero.
\\The standard basis is $G=\{f_1,f_2,f_3,f_4,f_5,f_6,f_7\}$. $G$ is a standard basis of $I$ in $\mathbb{Z}[x,y,z]_{\langle x,y,z \rangle}$.
\\Note that $f_1-zf_5=15x^2-15x^2yz^{10}=15x^2(1-yz^{10})$. Therefore $15x^2\in I$. Similarly one can see that $28y^3,28y^2z^5,35yz^5$ and $35xz^{15}\in I$. This implies that $\{15x^2,3x^2y+7yz^5,4xy^2-5xz^{10},28y^3,28y^2z^5,35yz^5,35xz^{15}\}$ is a standard basis of $I$.
\end{exmp}

\begin{thm}
$($cf. $[1]$ Theorem $4.5.9$ page $251$ $)$ Let $R$ be a PID, and $I$ be an ideal of $R[x_1, \ldots ,x_n]_>$. Assume that $\{f_1, \ldots ,f_n\}$ is a standard basis for $I$. Let $LT(f_i)=c_ix^{\alpha_i}$, for a saturated subset $J$ of $\{1,\ldots ,s\}$, let $c_J=gcd(c_j|j\in J)$ and write $c_J=\sum_{j\in J}a_jc_j$ $($any such representation will do$)$. Also, let $x^{\alpha _J} =lcm(x^{\alpha _j}|j\in J)$. Then the set
$$\{f_J=\sum_{j\in J}a_j\frac{x^{\alpha_J}}{x^{\alpha_j}}f_j|J \, is \, a \, saturated \, subset \, of \, \{1, \ldots ,s\}\}$$
is a strong standard basis for $I$. In particular, every non-zero ideal in $R[x_1, \ldots ,x_n]_>$ has a strong standard basis.
\end{thm}

\section{Standard Bases In The Formal Power Series Rings}

Finally we want to apply our results to compute standard bases in the formal power series ring $R[[x_1, \ldots ,x_n]]$ with coefficients in a ring $R$. Let $>$ be a local degree ordering, i.e, $>$ is a local ordering and $x^\alpha > x^\beta $ implies $deg(x^\alpha)\leq deg(x^\beta )$. A non-zero element $f\in R[[x_1, \ldots ,x_n]]$ can be written as $\sum_{v=0}^\infty a_vx^{\alpha_v}$, $a_v\in R$, $a_0\neq 0$ and $x^{\alpha_v} > x^{\alpha_{v+1}}$ for all $v$. As in Definition $2.5$, we define $LM(f),\,LE(f),\,LT(f),\,LC(f)$ and $tail(f)$. As in Definition $2.11$, we define a standard basis (respectively a strong standard basis) of an ideal $I\subseteq R[[x_1, \ldots ,x_n]]$.
\begin{prop}
Let $I< R[x_1, \ldots ,x_n]$ and $G$ is a standard basis (respectively a strong standard basis) of $I$ with respect to $>$, where $>$ is a local degree orderig. Then $G$ is a standard basis (respectively a strong standard basis) of $IR[[x_1, \ldots ,x_n]]$.
\end{prop}
\begin{proof}
Let $\{g_1,\ldots ,g_s\}$ be a standard basis of $I$ and $\overline{g}=\sum_{i=1}^s\overline{a_i}g_i \in IR[[x_1, \ldots ,x_n]]$, $\overline{g}\neq 0$. Let $c$ be an integer such that $LM(\overline{g})\notin \langle x_1,\ldots ,x_n\rangle^c$. Choose $a_i\in R[x_1, \ldots ,x_n]$ such that $\overline{a_i}-a_i\in \langle x_1,\ldots ,x_n\rangle^c$. Let $g=\sum_{i=1}^sa_ig_i$. Then $g\in I$ and $\overline{g}-g\in \langle x_1,\ldots ,x_n\rangle^c$. This implies $LT(\overline{g})=LT(g)$. If $G$ is a strong standard basis for $I$ then there exists $i$ such that $LT(g_i)|LT(g)=LT(\overline{g})$, i.e, $G$ is also a strong standard basis of $IR[[x_1, \ldots ,x_n]]$. Similarly it follows that a standard basis of $I$ is a standard basis of $IR[[x_1, \ldots ,x_n]]$.
\end{proof}

\section{Procedures}
Let $I< R[x_1, \ldots, x_n]$ be an ideal and $G=\{g_1,\ldots ,g_m\}$ be a standard basis of $I$. Then we can compute a strong standard basis of $I$ using the SINGULAR-procedures below.
\\ \\ \texttt{LIB"poly.lib";}
\\ \\ \texttt{proc powerSet(int n)} \indent \quad  //computes the set of all subsets of $\{1, \ldots, n\}$
\\ $\{$
\\ \indent \quad   \texttt{list L,K,S;}
\\ \indent \quad  \texttt{int i;}
\\  \indent \quad   \texttt{if(n==0)}
  \\  \indent \quad $\{$
    \\  \indent \quad  \indent \quad   \texttt{L[1]=L;}
     \\ \indent \quad  \indent \quad  \texttt{return(L);}
 \\ \indent \quad  $\}$
 \\ \indent \qquad  \texttt{if(n==1)}
 \\  \indent \qquad  $\{$
\\   \indent \quad  \indent \quad    \texttt{L[1]=L;}
  \\ \indent \quad  \indent \quad    \texttt{L[2]=list(1);}
    \\ \indent \quad  \indent \quad  \texttt{return(L);}
 \\ \indent \qquad $\}$
\\ \indent \qquad  \texttt{S=powerSet(n-1);}
\\ \indent \qquad  \texttt{int r=size(S);}
\\ \indent \qquad  \texttt{S[r+1]=list(n);}
\\ \indent \qquad  \texttt{for(i=2;i<=r;i++)}
\\ \indent \qquad  $\{$
\\ \indent \quad \indent \quad    \texttt{K=S[i];}
\\  \indent \quad \indent \quad   \texttt{K[size(K)+1]=n;}
\\  \indent \quad \indent \quad   \texttt{S[size(S)+1]=K;}
\\ \indent \qquad  $\}$
\\ \indent \qquad  \texttt{return(S);}
\\ $\}$
\\ \\ \texttt{proc satt(ideal I)} \indent \qquad   //computes the saturated subsets of $\{1, \ldots,size(I)\}$
\\$\{$ \indent \qquad   \indent \qquad \qquad \indent \qquad \indent   //w.r.t $\{LM(f)\,|\,f\in I\}$
\\  \indent \quad   \texttt{int j;}
\\ \indent \quad    \texttt{list K,J;}
\\ \indent \quad    \texttt{I=lead(I);}
\\ \indent \quad    \texttt{J=powerSet(size(I));}
\\ \indent \quad    \texttt{for(j=2;j<=size(J);j++)}
\\  \indent \quad   $\{$
\\  \indent \quad  \indent \quad     \texttt{ if(saturat(I,J[j]))}
\\  \indent \quad  \indent \quad     $\{$
\\   \indent \quad  \indent \quad  \indent \quad       \texttt{K[size(K)+1]=J[j];}
\\   \indent \quad  \indent \quad    $\}$
\\  \indent \quad    $\}$
\\   \indent \quad    \texttt{return(K);}
\\ $\}$
\\ \\ \texttt{proc specialGCD2(int a, int n)}
\\ $\{$
\\ \indent \quad   \texttt{int x=a mod n;}
\\  \indent \quad   \texttt{if(x==0){return(list(0,1,n));}}
\\  \indent \quad   \texttt{list L=specialGCD2(n,x);}
\\  \indent \quad   \texttt{return(list(L[2],L[1]-(a-x)*L[2]/n,L[3]));}
\\ $\}$
\\ \\ \texttt{proc specialGCD(list L)}  \indent \quad  //computes the gcd over $\mathbb{Z}$
\\ $\{$
\\  \indent \quad   \texttt{int i;}
\\  \indent \quad   \texttt{for(i=1;i<=size(L);i++){L[i]=int(L[i]);}}
\\  \indent \quad   \texttt{if(size(L)==1){return(L);}}
 \\ \indent \quad   \texttt{if(size(L)==2){return(specialGCD2(L[1],L[2]));}}
 \\ \indent \quad   \texttt{bigint p=L[size(L)];}
\\  \indent \quad   \texttt{L=delete(L,size(L));}
\\  \indent \quad  \texttt{list T=specialGCD(L);}
\\ \indent \quad    \texttt{list S=specialGCD2(T[size(T)],p);}
\\  \indent \quad   \texttt{for(i=1;i<=size(T)-1;i++)}
\\ \indent \quad    $\{$
\\ \indent \quad   \indent \quad \texttt{ T[i]=T[i]*S[1];}
\\  \indent \quad   $\}$
 \\  \indent \quad  \texttt{p=T[size(T)];}
\\   \indent \quad   \texttt{T[size(T)]=S[2];}
\\   \indent \quad   \texttt{T[size(T)+1]=S[3];}
\\  \indent \quad   \texttt{return(T);}
\\ $\}$
\\ \\ \texttt{proc coeffJ(list C,list J)}  \indent \quad // computes the coefficients $a_j$ such that \\ $\{$\indent \qquad  \indent \qquad  \indent \qquad \indent \qquad \indent \qquad // $gcd(c_j|j\in J)=\sum_{j\in J}a_jc_j$
\\  \indent \quad    \texttt{if(size(J)==1){return(list(1));}}
\\  \indent \quad    \texttt{int n=size(J);}
\\  \indent \quad    \texttt{int i;}
\\ \indent \quad     \texttt{list L,M;}
\\  \indent \quad    \texttt{for(i=1;i<=n;i++)}
\\  \indent \quad    $\{$
\\  \indent \quad   \indent \quad \texttt{L[size(L)+1]=C[J[i]];}
\\  \indent \quad   $\}$
\\  \indent \quad   \texttt{L=specialGCD(L);}
\\  \indent \quad   \texttt{L=delete(L,size(L));}
\\   \indent \quad   \texttt{return(L);}
\\ $\}$
\\ \\ \texttt{proc lcmJ(ideal X,list J)} \indent \quad  \indent \quad   // computes $lcm(f_j\,|\,f_j\in X,\, j\in J)$
\\$\{$
\\  \indent \quad \texttt{poly p=X[J[1]];}
\\   \indent \quad \texttt{int i;}
\\  \indent \quad  \texttt{for(i=2;i<=size(J);i++)}
\\  \indent \quad $\{$
\\  \indent \quad \indent \quad \texttt{ p=lcmS(p,X[J[i]]);}
\\  \indent \quad $\}$
\\  \indent \quad \texttt{return(p);}
\\ $\}$
\\ \\ \texttt{proc maxZ(intvec a,intvec b)} \indent \qquad  // computes the maximum of two
\\$\{$  \indent \qquad  \indent \qquad \indent \qquad  \indent \qquad \qquad //integer vectors
\\  \indent \quad  \texttt{int i,j;}
\\  \indent \quad \texttt{intvec c;}
\\  \indent \quad \texttt{for(i=1;i<=size(a);i++)}
\\  \indent \quad $\{$
 \\  \indent \quad \indent \quad \texttt{if(a[i]>=b[i])}
\\    \indent \quad \indent \quad $\{$
\\   \indent \quad   \indent \quad  \indent \quad \texttt{c[i]=a[i];}
\\  \indent \quad  \indent \quad   $\}$
\\   \indent \quad    \indent \quad       \texttt{else}
\\   \indent \quad   \indent \quad         $\{$
\\   \indent \quad   \indent \quad   \indent \quad \texttt{c[i]=b[i];}
\\     \indent \quad   \indent \quad $\}$
\\   \indent \quad   $\}$
\\  \indent \quad  \texttt{ return(c);}
\\ $\}$
\\ \\ \texttt{proc lcmS(poly p,poly q)}\indent \qquad  \qquad // computes the LCM of two monomials
\\  \indent \quad    \texttt{intvec a=leadexp(p);}
\\  \indent \quad    \texttt{intvec b=leadexp(q);}
\\  \indent \quad    \texttt{intvec c=maxZ(a,b);}
\\ \indent \quad    \texttt{ int i,j;}
\\ \indent \quad    \texttt{ poly s=1;}
\\ \indent \quad     \texttt{for(i=1;i<=nvars(basering);i++)}
\\  \indent \quad   $\{$
\\  \indent \quad   \indent \quad  \texttt{s=s*var(i)}\verb"^"\texttt{c[i]};
\\  \indent \quad    $\}$
\\   \indent \quad    \texttt{ return(s);}
\\ $\}$
\\ \\ \texttt{proc leadTerm(ideal I)}\indent \qquad // give a list containing two lists, first \\ $\{$ \indent \qquad  \qquad \indent \qquad  \qquad \indent \quad //contains leading coefficients and second contains
\\  \indent \quad    \texttt{int i;} \indent \qquad  \qquad \indent \qquad  \qquad //leading monomials
\\ \indent \quad     \texttt{list L;}
\\ \indent \quad     \texttt{ideal J;}
\\  \indent \quad    \texttt{I=lead(I);}
\\  \indent \quad    \texttt{for(i=1;i<=size(I);i++)}
\\  \indent \quad    $\{$
\\  \indent \quad \indent \quad   \texttt{L[size(L)+1]=leadcoef(I[i]);}
\\  \indent \quad  \indent \quad   \texttt{J[size(J)+1]=leadmonom(I[i]);}
\\  \indent \quad   $\}$
\\  \indent \quad     \texttt{ return(list(L,J));}
\\ $\}$
\\ \\ \texttt{proc clean(ideal G)} \indent \qquad// delete from an ideal $f_i$ such that \\ $\{$ \indent \qquad  \indent \qquad \indent \qquad // there exists $f_j\in I$ with $LM(f_j)$ divides $LM(f_i)$
\\ \indent \quad   \texttt{int i,j;}
\\  \indent \quad    \texttt{while(i<=size(G)-1)}
\\  \indent \quad   $\{$
\\  \indent \quad   \indent \quad \texttt{i++;}
\\  \indent \quad   \indent \quad \texttt{j=i;}
 \\  \indent \quad   \indent \quad  \texttt{while(j} \verb"<" \texttt{size(G))}
\\  \indent \quad   \indent \quad    $\{$
\\  \indent \qquad   \indent \qquad   \texttt{ j++;}
 \\  \indent \qquad  \indent \qquad   \texttt{if((G[i]!=0)}\verb"&&"\texttt{(G[j]!=0))}
\\  \indent \qquad  \indent \qquad   $\{$
 \\  \indent \qquad  \indent \qquad \qquad   \texttt{if (lead(G[j])/lead(G[i])!=0)}
 \\  \indent \qquad  \indent \qquad \qquad  $\{$
\\ \indent \qquad   \indent \qquad  \qquad   \indent \qquad \texttt{G[j]=0;}
 \\ \indent \qquad   \indent \qquad   \qquad  \indent \qquad \texttt{G=simplify(G,2);}
 \\  \indent \qquad   \indent \qquad   \qquad  \indent \qquad \texttt{j--;}
\\ \indent \qquad   \indent \qquad \qquad  $\}$
\\  \indent \qquad  \indent \qquad \indent \qquad \texttt{if(G[j]!=0)}
\\  \indent \qquad  \indent \qquad \indent \qquad $\{$
 \\  \indent \qquad  \indent \qquad \indent \qquad  \qquad \texttt{if (lead(G[i])/lead(G[j])!=0)}
 \\  \indent \qquad  \indent \qquad   \indent \qquad  \qquad  $\{$
\\  \indent \qquad  \indent \qquad  \indent \qquad  \indent \qquad  \qquad  \texttt{G[i]=0;}
\\  \indent \qquad  \indent \qquad  \indent \qquad  \indent \qquad  \qquad \texttt{G=simplify(G,2);}
 \\  \indent \qquad   \indent \qquad  \indent \qquad  \indent \qquad \qquad  \texttt{i--;}
 \\  \indent \qquad  \indent \qquad  \indent \qquad  \indent \qquad  \qquad  \texttt{break;}
 \\  \indent \qquad   \indent \qquad   \indent \qquad \qquad  $\}$
\\  \indent \qquad  \indent \qquad \indent \qquad $\}$
 \\  \indent \qquad  \indent \qquad    $\}$
\\  \indent \quad   \indent \quad  $\}$
\\  \indent \quad    $\}$
\\  \indent \quad   \indent \quad \texttt{return(simplify(G,2));}
\\ $\}$
\\ \\ \texttt{proc strongSB(ideal I)} \indent \qquad // computes a strong standard basis
\\ $\{$ \indent \qquad \indent \qquad \indent \qquad \qquad // for an ideal if standard basis is given
\\ \indent \quad   \texttt{def R=basering;}
\\  \indent \quad   \texttt{list rl=ringlist(R);}
\\  \indent \quad   \texttt{rl[1]=0;}
\\  \indent \quad   \texttt{def S=ring(rl);}
\\ \indent \quad    \texttt{setring S;}
\\ \indent \quad    \texttt{ideal I=imap(R,I);}
\\  \indent \quad   \texttt{list L=leadTerm(I);}
\\  \indent \quad   \texttt{list C=L[1];}
\\  \indent \quad   \texttt{ideal X=L[2];}
\\  \indent \quad   \texttt{list J=satt(X);}
\\ \indent \quad    \texttt{int i,j;}
\\ \indent \quad   \texttt{ideal G;}
\\ \indent \quad    \texttt{list P;}
\\  \indent \quad   \texttt{poly f,q;}
\\ \indent \quad   \texttt{for(i=1;i<=size(J);i++)}
\\ \indent \quad    $\{$
\\  \indent \quad   \indent \quad \texttt{P=coeffJ(C,J[i]);}
\\  \indent \quad   \indent \quad  \texttt{q=lcmJ(X,J[i]);}
\\   \indent \quad   \indent \quad \texttt{f=0;}
\\  \indent \quad   \indent \quad \texttt{for(j=1;j<=size(J[i]);j++)}
\\  \indent \quad   \indent \quad  $\{$
\\  \indent \quad   \indent \quad   \indent \quad \texttt{f=f+q/X[J[i][j]]*I[J[i][j]]*P[j];}
\\  \indent \quad   \indent \quad     $\}$
\\    \indent \quad   \indent \quad\texttt{G[size(G)+1]=f;}
\\  \indent \quad   $\}$
\\  \indent \quad    \texttt{setring R;}
\\  \indent \quad    \texttt{ideal G=imap(S,G);}
\\ \indent \quad     \texttt{return(clean(G));}
\\ $\}$
\begin{exmp}
Consider Example $6.3$, we now compute strong standard basis for $I=\{15x^2,3x^2y+7yz^5,4xy^2-5xz^{10},28y^3,28y^2z^5,35yz^5,35xz^{15}\}$.
\\ \verb">" $ring\, R=integer,(x,y,z),ds$;
\\ \verb">" $ideal\, I=15x2,3x2y+7yz5,4xy2-5xz10,28y3,28y2z5,35yz5,35xz15$;
\\ \verb">" $strongSB(I)$;
\\ $[1]=15x2$
\\ $[2]=3x2y+7yz5$
\\ $[3]=4xy2-5xz10$
\\ $[4]=x2y2-7y2z5-5x2z10$
\\ $[5]=28y3$
\\ $[6]=35yz5$
\\ $[7]=x2yz5+84yz10$
\\ $[8]=7y2z5$
\\ $[9]=35xz15$
\\ $[10]=5x2z15$



\end{exmp}

\end{document}